\documentclass[10pt,twocolumn]{amsart}
\usepackage{etex} 
\usepackage{amssymb}
\usepackage[noadjust]{cite}
\usepackage{booktabs}
\usepackage{xypic}
\usepackage{url}
\usepackage{hyphenat}
\usepackage{mathtools}
\usepackage{cancel} 
\usepackage[T1]{fontenc}
\usepackage[utf8]{inputenc}
\usepackage{mdwtab}
\usepackage{enumitem}
\usepackage{bbding}
\usepackage{dsfont}
\usepackage{hyphenat}
\DeclareMathAlphabet\mathbfcal{OMS}{cmsy}{b}{n}

\usepackage{mathrsfs}
\usepackage{xr-hyper}
\usepackage{xcolor}

\definecolor{darkblue}{rgb}{0,0,0.4}
\usepackage[pdftex,lmargin=1in, rmargin=1in, tmargin=1in, bmargin=1in]{geometry}
  \usepackage[bookmarks=true, bookmarksopen=true,%
    bookmarksdepth=3,bookmarksopenlevel=2,%
    colorlinks=true,%
    linkcolor=darkblue,%
    citecolor=darkblue,%
    filecolor=darkblue,%
    menucolor=darkblue,%
    urlcolor=darkblue]{hyperref}
  \hypersetup{pdftitle={Floer homology beyond borders}}
  \hypersetup{pdfauthor={Robert Lipshitz, Peter S. Ozsváth and Dylan
      P. Thurston}}
\usepackage{color}
\usepackage{tikz}

\usetikzlibrary{matrix,arrows,cd}
\tikzstyle{every picture}=[> = latex']

\tikzset{taa/.style={->, double}}
\tikzset{moda/.style={->, dashed}}
\tikzset{alga/.style={->, thick}}
\newread\testin

\graphicspath{{draws/}{mpdraws/}{}}
\makeatletter
\def\input@path{{}{draws/}}
\makeatother

\makeatletter
\newcommand\mi@kern[1]{%
  \settowidth\@tempdima{$\mi@obj^{#1}$}
  \kern-\@tempdima
  #1
  \settowidth\@tempdima{$\mi@obj$}
  \kern\@tempdima
}

\newtoks\mi@toksp
\newtoks\mi@toksb
\DeclareRobustCommand{\manyindices}[5]{
  \def\mi@obj{#5}
  \mi@toksp\expandafter{\mi@kern{#2}}
  \mi@toksb\expandafter{\mi@kern{#1}}
  \@mathmeasure4\textstyle{#5_{#1}^{#2}}
  \@mathmeasure6\textstyle{#5_{#3}^{#4}}
  \dimen0-\wd6 \advance\dimen0\wd4
  \@mathmeasure8\textstyle{\hphantom{{}_{#1}^{#2}}#5^{\the\mi@toksp#4}_{\the\mi@toksb#3}}
  \hbox to \dimen0{}{\kern-\dimen0\box8}
}
\makeatother 

\usepackage{ifpdf}
\ifpdf
  
\else
  
\fi

\newcommand{\RR}{\mathbb R}

\newcommand{\ZZ}{\mathbb Z}

\newcommand{\FF}{\mathbb F}

\newcommand{\co}{\nobreak\mskip2mu\mathpunct{}\nonscript
  \mkern-\thinmuskip{:}\penalty300\mskip6muplus1mu\relax}

\newcommand{\bdy}{\partial}

\newcommand{\lbracket}{[}
\newcommand{\rbracket}{]}

\theoremstyle{plain}

\numberwithin{equation}{section}

\theoremstyle{definition}

\theoremstyle{remark}

\newcommand{\HF}{\mathit{HF}}
\newcommand{\HFa}{\widehat {\HF}}

\newcommand{\HFm}{{\HF}^-}

\newcommand{\CFa}{\widehat {\mathit{CF}}}
\newcommand{\CFm}{\mathit{CF}^-}

\newcommand{\HFKa}{\widehat{\mathit{HFK}}}

\newcommand\HH{\mathit{HH}}

\newcommand\Hochschild\HH

\newcommand{\Ainf}{A_\infty}
\newcommand{\Alg}{\mathcal{A}}

\newcommand{\DD}{\textit{DD}}

\newcommand{\CFD}{\mathit{CFD}}
\newcommand{\CFDD}{\mathit{CFDD}}
\newcommand{\CFA}{\mathit{CFA}}

\newcommand{\CFDA}{\mathit{CFDA}}
\newcommand{\CFDAa}{\widehat{\CFDA}}

\newcommand{\CFDa}{\widehat{\CFD}}

\newcommand{\CFDm}{\CFD^-}

\newcommand{\HFK}{\mathit{HFK}}
\newcommand{\HFKm}{\HFK^-}
\newcommand{\CFDDa}{\widehat{\CFDD}}
\newcommand{\CFDDm}{\CFDD^-}

\newcommand{\CFAa}{\widehat{\CFA}}

\newcommand{\dg}{\textit{dg} }

\newcommand\Id{\mathbb{I}}

\newcommand\DT{\boxtimes}

\newcommand\Zmod[1]{\mathbb{Z}/{#1}\mathbb{Z}}

\DeclareMathOperator{\nbd}{nbd}

\DeclareMathOperator{\Mor}{Mor}

\makeatletter
\newcommand\honestalg[3]{\bigl\lbracket
\begin{smallmatrix} #1\@ifempty{#3}{}{&#3} \\ #2 \end{smallmatrix}
\bigr\rbracket}

\makeatother

\newcommand{\lsup}[2]{{}^{#1}\mskip-.6\thinmuskip#2}

\newcommand{\MAlg}{\Alg_-}

\newcommand{\CFAm}{\CFA^-}

\newcommand{\kotimes}[1]{\otimes}

\newcommand\goesto\mapsto

\definecolor{darkgreen}{rgb}{0,0.4,0} 
\definecolor{darkbrown}{rgb}{.48,0.33,.24}  

\begin{document}
\title[Floer homology beyond borders]{Floer homology beyond borders}
\author[Lipshitz]{Robert Lipshitz}
\thanks{RL was partly supported by NSF grant DMS-2204214 and by a Simons Fellowship.}
\address{Department of Mathematics, University of Oregon\\
  Eugene, OR 97403}
\email{lipshitz@uoregon.edu}

\author[Ozsv\'ath]{Peter Ozsv\'ath}
\thanks{PSO was supported by NSF grant DMS-2104536.}
\address {Department of Mathematics, Princeton University\\ New
  Jersey, 08544}
\email {petero@math.princeton.edu}

\author[Thurston]{Dylan~P.~Thurston}
\thanks {DPT was supported by NSF grant DMS-2110143.}
\address{Department of Mathematics\\
         Indiana University,\linebreak
         Bloomington, Indiana 47405\\
         USA}
\email{dpthurst@indiana.edu}

\begin{abstract}
  Bordered Floer homology is an invariant for 3\hyp manifolds with
  boundary, defined by the authors in 2008.
  It extends the Heegaard
  Floer homology of closed 3\hyp manifolds, defined in earlier work of
  Zolt{\'a}n Szab{\'o} and the second author. In addition to its
  conceptual interest, bordered Floer homology also provides powerful
  computational tools. This survey outlines the
  theory, focusing on recent developments and applications.
\end{abstract} 

\maketitle



The goal of this note is to survey some developments in bordered Floer
homology, an extension of Heegaard Floer homology. We start with a
biased outline of Heegaard Floer homology, focusing somewhat on the
aspects relevant to the extension. We then briefly outline the
structure of bordered Floer homology, before turning to a discussion
of its use for computations, extensions of it, and some recent
applications to 4-dimensional topology and its 3-dimen\-sional
shadows.

\section{Heegaard Floer homology}\label{sec:Heegaard-Floer}

Heegaard Floer homology~\cite{OS04:HolomorphicDisks,OS06:HolDiskFour}
is an invariant of 3- and 4\hyp dimensional manifolds
defined by Zolt{\'a}n
Szab{\'o} and the second author via methods from
symplectic geometry (specifically, Lagrangian Floer homology~\cite
{Floer88:LagrangianHF}). The
construction was inspired by gauge theory, especially the
Seiberg-Witten invariants~\cite{Witten}, which also give a package of
invariants with similar properties~\cite{KronheimerMrowka}.
Thanks to work of {\c{C}}a\u{g}atay Kutluhan, Yi-Jen Lee, and Clifford
H. Taubes~\cite{KutluhanLeeTaubes20:HFHMI,KutluhanLeeTaubes20:HFHMII,KutluhanLeeTaubes20:HFHMIII,KutluhanLeeTaubes20:HFHMIV,KutluhanLeeTaubes20:HFHMV}
and Vincent Colin, Paolo Ghiggini, and Ko
Honda~\cite{ColinGhigginiHonda11:HF-ECH-1,ColinGhigginiHonda11:HF-ECH-2,ColinGhigginiHonda11:HF-ECH-3},
building on earlier work of Michael Hutchings and
Taubes~\cite{Hutchings02:ECH,HutchingsTaubes07:GluingI,HutchingsTaubes09:GluingII},
we now know that, at least for 3\hyp manifolds, these two differently
defined invariants agree.  The two
perspectives have different relative strengths: Seiberg-Witten theory
is more directly connected with the differential geometry of the
underlying manifold, whereas Heegaard Floer homology is more closely
connected to its combinatorial topological properties.

On a formal level, Heegaard Floer homology associates to a 3\hyp manifold
$Y$ the chain homotopy type of a chain complex $\CFm(Y)$ over a
polynomial algebra
in an indeterminant $U$; in particular, its homology $\HFm(Y)$ is also a 3\hyp manifold invariant.
The $U=0$
specialization of $\CFm(Y)$, denoted $\CFa(Y)$, and its homology,
$\HFa(Y)$, are also useful invariants.
(The invariant $\HFm(Y)$
is analogous to the $S^1$-equivariant homology of a space,
while $\HFa(Y)$ is like its non-equivariant homology; cf.~\cite{Donaldson96:survey,Manolescu03:SW-spectrum-1,ManolescuKronheimerSWspectrum,KronheimerMrowka}.)

A striking
result of Sucharit Sarkar and Jiajun Wang led to a combinatorial
scheme for computing $\HFa(Y)$, via so-called \emph{nice} auxiliary
data~\cite{SarkarWang07:ComputingHFhat}, a notable achievement given
the many 3\hyp dimensional applications of $\HFa(Y)$.
However, the construction of invariants of smooth, closed
4\hyp dimensional manifolds requires the use the $U$-unspecialized
version~$\HFm(Y)$ \cite{OS06:HolDiskFour}. Specifically, a cobordism
between two
3\hyp manifolds $Y_1$ and $Y_2$ gives rise to a map between between
their respective Heegaard Floer homology groups. The closed, smooth
4\hyp manifold invariant is constructed using the interaction between
these maps with the $U$-module structure of $\HFm(Y)$. Although the
4\hyp manifold invariant can now be, in principle, described
combinatorially~\cite{MOT:grid}, the algorithm for computing it is still
too unwieldly to be implemented in practice.

\subsection{Knot invariants from Heegaard Floer homology}

In 2003, Jacob Rasmussen and, independently, Szab{\'o} and the second
author constructed a version of Heegaard Floer homology for knots in
$S^3$. The resulting {\em knot Floer homology} is a bigraded abelian group
whose graded Euler characteristic is the Alexander polynomial. This
invariant contains much topological information about the underlying knot,
including its Seifert genus~\cite{OS04:ThurstonNorm} and whether its
complement fibers over the circle~\cite{Ghiggini08:FiberedGenusOne,Ni07:FiberedKnot};
see also~\cite{Juhasz08:SuturedDecomp}. 

Building on Sarkar's earlier work, Ciprian Manolescu,
Sarkar, and
the second author gave a combinatorial construction of knot
Floer homology, now called {\em grid
  homology}~\cite{MOS06:CombinatorialDescrip,MOST07:CombinatorialLink,OSS15:grid-book}. For a knot with $n$ crossings, grid
homology can be computed as the homology groups of a chain complex
with roughly $n!$ generators. As such, it can be computed explicitly
for knots with $\leq 16$ or so
crossings~\cite{BaldwinGillam12:compute,Beliakova10:grid}.

Heegaard Floer
homology can also be used to construct link invariants in a different
manner.  Given a link $L\subset S^3$ there is a closed 3\hyp manifold,
the \emph{double cover of $S^3$ branched along $L$}, which is
the unique 3\hyp manifold $\Sigma(L)$ that admits a map $f\co
\Sigma(L)\to S^3$ which is a $2$-to-$1$ covering space away from $L$
and is modeled on $z\mapsto z^2$ in the normal
direction to $L$.  Thus, given a closed 3\hyp manifold
invariant, applying that invariant to $\Sigma(L)$ gives a link
invariant. In the case of $\HFa$, the link invariant
$\HFa(\Sigma(L))$ is surprisingly similar to an
invariant defined by Mikhail Khovanov, inspired by the Jones
polynomial and constructions in
representation theory~\cite{Khovanov00:CatJones}.
For an unlink, the (reduced) Khovanov homology and the Heegaard Floer
homology of its branched double cover agree. More strikingly, both
$\HFa(\Sigma(L))$ and the Khovanov homology of $L$ satisfy an exact
triangle when one resolves a crossing in the two possible ways; for
Heegaard Floer homology, this is a special case of the \emph{surgery
  exact triangle}~\cite{OS04:HolDiskProperties},
parallel to one envisaged by Andreas Floer in gauge
theory~\cite{FloerTriangles,BraamDonaldson95:triangle}. These
properties lead to a precise relationship between $\HFa(\Sigma(L))$
and Khovanov homology: for every link diagram, there is a spectral
sequence whose $E_2$ term is the
(reduced) Khovanov homology of $L$ and which converges to
$\HFa(\Sigma(L))$ (both with mod-$2$
coefficients)~\cite{BrDCov}. 
By now we know that there are a number of other spectral sequences from Khovanov
homology to various Floer-homological invariants: sutured
Heegaard Floer homology~\cite{GrigsbyWehrli:detects}, instanton Floer homology~\cite{KronheimerMrowka11:detect}, monopole Floer
homology~\cite{Bloom11:ss}, and knot Floer
homology~\cite{Dowlin:kh-to-hfk}. Because of the topological content of these latter
invariants, the spectral sequences have interesting consequences for Khovanov homology.

\subsection{Contact geometry and Heegaard Floer homology}

Gauge theory is closely tied with the symplectic geometry of the underlying manifold. According to a
celebrated theorem of Donaldson's~\cite{DonaldsonPolynomials}, his
smooth 4\hyp manifold invariant is non-zero for a K{\"a}hler
manifold. Taubes~\cite{Taubes94:SW-symplectic} proved the analogous
non-vanishing theorem for Seiberg-Witten invariants of symplectic
manifolds. This has a 3\hyp dimensional counterpart: Kronheimer and
Mrowka~\cite{KM97:contact,KMOS} constructed an invariant for contact
structures over 3\hyp manifolds, which takes values in its
Seiberg-Witten Floer homology; moreover, this contact invariant is
non-trivial when the contact structure is fillable, in a suitable
sense, by a symplectic manifold. These results,
combined with work of Yakov Eliashberg and William
Thurston~\cite{ET98:confoliations} relating contact structures and
foliations, form a bridge
between gauge theory and fundamental 
topological constructions developed by David Gabai~\cite{Gabai83:foliations,Gabai84:knot-foliations}.
This bridge was exploited by Kronheimer and Mrowka in their celebrated proof
that all knots in $S^3$ have ``Property P''~\cite{KM04:P}.

An analogous bridge between Heegaard Floer homology and contact
geometry is built upon the work of Emannuel
Giroux~\cite{Giroux02:correspondence-ICM}, who reformulated contact
structures over $Y$ as certain equivalence classes of open book
decompositions of $Y$. Using Giroux's correspondence, Szab{\'o} and
the second author were able to define an invariant for contact
structures over $Y$, analogous to the
Kronheimer-Mrowka contact invariant, but taking values in
$\HFa(-Y)$~\cite{OS05:Contact} (where $-Y$ denotes orientation
reverse). By work of Taubes \cite{Taubes10:ECH-SW1,Taubes10:ECH-SW2,Taubes10:ECH-SW3,Taubes10:ECH-SW4,Taubes10:ECH-SW5} and
Colin-Ghiggini-Honda~\cite{ColinGhigginiHonda11:HF-ECH-1,ColinGhigginiHonda11:HF-ECH-2,ColinGhigginiHonda11:HF-ECH-3},
this contact invariant is a refinement of Kronheimer and Mrowka's. (See
also~\cite{BS21:same-contact}.)

This contact invariant forms an important technical tool within the
theory
(used, e.g., to prove that $\HFKa$ detects the Seifert
genus~\cite{OS04:ThurstonNorm}); and it is also a vehicle for studying
contact
structures in their own right (e.g.,~\cite{LScontactInv0,LScontactInvI,LScontactInvII}).

\section{Bordered Floer homology}

\emph{Bordered Floer homology} is an invariant for 3\hyp manifolds with
boundary~\cite{LOT1}. The basic structure is easy to describe.
To a surface $F$, it associates a differential
graded algebra $A(F)$. Given an $F$-bordered 3\hyp manifold, that is, an oriented
3\hyp manifold $Y$ and an orientation-preserving
diffeomorphism $\phi\co F\to \partial Y$, the theory associates a
right $\Ainf$-module $\CFAa(Y,\phi)_{A(F)}$ over the algebra
$A(F)$ and
a left \dg module of a particular kind (a twisted complex or, in the internal
language of the subject, \emph{type $D$ structure})
$\lsup{A(-F)}\CFDa(Y,\phi)$ over the algebra associated to the
orientation-reversed surface. Both $\CFAa$ and $\CFDa$ depend on some
auxiliary choices, but up to
homotopy equivalence are independent of those choices.
Given two bordered 3\hyp manifolds $(Y_1,\phi_1\co -F\to \partial Y_1)$ and
$(Y_2,\phi_2\co F\to \partial Y_2)$, we can recover the Heegaard Floer
invariant $\CFa(Y_1\cup_{\phi_2\circ\phi_1^{-1}} Y_2)$ as either a tensor product
\[
  \CFAa(Y_1)_{A(F)}\DT\lsup{A(F)} \CFDa(Y_2)
\]
or as the complex of $\Ainf$-module morphisms
\[
  \Mor^{A(F)}(\CFAa(-Y_1),\CFAa(Y_2))  
\]
or type $D$ structure morphisms
\[
  \Mor_{A(-F)}(\CFDa(-Y_1),\CFDa(Y_2))
\]
\cite{LOT1,LOTHomPair,Auroux10:Bordered}. We call these results
\emph{pairing theorems} (to distinguish them from the gluing theorems
used in the analysis underlying the subject).

The analytic proof of the pairing theorem involves a neck-stretching
argument, inspired by the proofs of product formulas for
gauge-theoretic
invariants~\cite{DonaldsonFloer,KronheimerMrowka}. (See
also~\cite{EGH00:IntroductionSFT,BEHWZ03:CompactnessInSFT}.) We start
by taking a Heegaard diagram which is adapted to the splitting
$Y_1\cup_{F} Y_2$ and degenerating the Heegaard surface along a
corresponding
circle. The counts of holomorphic disks that constitute the Heegaard
Floer differential degenerate into counts of pairs of holomorphic
disks on the two sides, satisfying
a matching condition along $F$. So far, the
counts are not algebraic: because of the matching condition,
high-dimensional moduli spaces on the two sides contribute to the
count. The next step is to deform the matching condition. Formally,
this deformation can be realized as a cellular approximation to the
diagonal in James Stasheff's
associahedra~\cite{Stasheff63:associahedron1}.  After the deformation,
we obtain algebraic counts that correspond to the tensor product
description from the pairing theorem.

The dependence of $\CFAa(Y_1,\psi)$ on the choice of
boundary parameterization is governed by certain
bimodules~\cite{LOT2}.  Given a
surface diffeomorphism $\phi\co
F\to F$, there is a bimodule
$\CFDAa(\phi)$, with the property that
for any bordered 3\hyp manifold $(Y_1,\psi\co F\to Y_1)$
\begin{equation}
  \label{eq:DependsOnParameterization}
  \CFAa(Y_1,\psi)\DT \CFDAa(\phi)
  \simeq \CFAa(Y_1,\psi\circ\phi).
\end{equation}
Moreover, the tensor product of the bimodules associated to $\phi_1$ and $\phi_2$ coincides with the bimodule associated to the composite $\phi_2\circ \phi_1$.

The algebra $A(F)$ is defined combinatorially, from a handle
decomposition for $F$. (Different handle decompositions lead to
derived equivalent, but non-quasi-isomorphic, algebras~\cite{LOT2}.)
By contrast, $\CFAa(Y)$ and $\CFDa(Y)$ are defined by counting
$J$-holomorphic curves. For suitable
(\emph{nice}, extending~\cite{SarkarWang07:ComputingHFhat}) choices of
auxiliary data, these curve counts are combinatorial, and there is
also a combinatorial proof of the tensor-product pairing
theorem. However, there is no known combinatorial proof of invariance
for $\CFAa$ and $\CFDa$ along these lines (but
see~\cite{Zhan16:proofs}), and nice auxiliary choices result in
enormous complexes. So, in practice---and for most of the
computational applications below---one often ends up counting
$J$-holo\-morphic curves, not using nice auxiliary choices.

\section{Computations from bordered Floer homology}

Every closed 3\hyp manifold $Y^3$ admits a \emph{Heegaard decomposition}
$Y=H_1\cup_{F} H_2$ into two handlebodies glued together along their
boundaries. If we identify each $H_i$ with some standard (bordered)
handlebody, then $Y$ is determined by the gluing diffeomorphism, which
is a map $\phi\co F\to F$. For an appropriate choice of standard
handlebody, it is easy to compute the bordered invariants $\CFDa(H_i)$
and $\CFAa(H_i)$. So, in view of
Equation~\eqref{eq:DependsOnParameterization}, to compute $\HFa(Y)$
for general $Y$, it suffices to compute the bimodules associated to
surface diffeomorphisms---or, in fact, for any set of generators of
the mapping class group. It turns out that there are particularly
simple generators, called \emph{arcslides}, for a groupoid extension of the
mapping class group~\cite{AndersenBenePenner09:MCGroupoid,Bene08:ChordDiagrams}, for which the type \DD\ bimodules are determined
by a few simple curve counts and the structure equation
$\bdy^2=0$~\cite{LOT4}.
A computer implementation of this algorithm~\cite{Zhan:code} (and
refinements of it~\cite{Zhan14:thesis}) is practical for many manifolds.

The bordered description arising from factoring mapping classes
contains information beyond $\HFa$ of the underlying 3\hyp manifold. When
$Y=\Sigma(L)$ is a branched double cover of a link $L$ in $S^3$, we
can factor its gluing map as a product of Dehn twists along an
explicit set of curves which correspond to crossings in a projection
of $L$. The bimodules for these Dehn twists can be expressed as
mapping cones between pairs of bimodules, corresponding to the two
resolutions of the crossing. This description induces a filtration on
$\CFa(\Sigma(L))$~\cite{LOT:DCov1}, which can be
shown~\cite{LOT:DCov2} to induce the aforementioned spectral
sequence~\cite{BrDCov} from the Khovanov homology of $L$ to
$\HFa(\Sigma(L))$. Like the arcslide bimodules,
these resolution bimodules and the maps between them can be
described explicitly, using only a few simple curve counts and the
structure equation; this gives a combinatorial formula for the
spectral sequence \cite{LOT:DCov1}. (Another conjectural combinatorial
description of the spectral sequence was given by
Szab\'o~\cite{Szabo15:ss}; it would be interesting to relate the two.)
The higher differentials in the original spectral sequence count
holomorphic polygons, and the key step in showing the spectral
sequences agree is an extension of the pairing theorem to
polygons~\cite{LOT:DCov2}.

\subsection{Satellite knots}
Suppose $K$ is a knot in $S^3$, say, and $P$ is another knot, called
the {\em pattern}, in the solid torus $S^1\times D^2$.  Choose a
framing for $K$ or, equivalently, an identification of
$\bdy(S^3\setminus \nbd(K))$ with $\bdy(S^1\times D^2)$. Replacing
$\nbd(K)$ with the copy of $S^1\times D^2$ containing~$P$ using this framing gives a
new knot $K_P$ in $S^3$, called a {\em satellite knot}.

It is natural to wonder how invariants of $K$, $P$, and $K_P$
relate. For the Alexander polynomial, there is a simple formula:
\[
\Delta_{K_P}(t)=\Delta_K(t^n)\Delta_P(t)
\]
where $\Delta_P$ is the Alexander polynomial of $P\subset S^3$,
and $n$ is the winding number of the pattern around the solid torus. 
For knot Floer homology, an answer comes from bordered Floer
homology. First, one can express the bordered invariants of the
exterior of $K$ in terms of the knot Floer homology of
$K$~\cite{LOT1} (see also~\cite{HeddenLevine16:splicing,Hockenhull:bord-link}). Then, the knot Floer homology of $K_P$ is a tensor
product of this bordered invariant with an appropriate module. In
fact, Ina Petkova showed that, in a precise sense, this construction
lifts the classical formula for the Alexander
polynomial~\cite{Petkova18:decat}. Satellite formulas from bordered
Floer homology have been used by many authors for computations and
applications; some applications of these will be mentioned below. We
also note that there is important earlier work on satellite operations
in knot Floer homology by Matthew
Hedden~\cite{Hedden,HeddenWhitehead} and Eaman
Eftekhary~\cite{Eftekhary05:LongitudeWhitehead}.

\section{Extensions of bordered Floer homology}
\subsection{Bordered-sutured Floer homology}
Inspired by Gabai's notion of sutured manifold
hierarchies~\cite{Gabai83:foliations,Gabai84:knot-foliations},
Andr\'as Juh\'asz introduced a different
extension \cite{Juhasz06:Sutured} of Heegaard Floer homology to
3\hyp manifolds with sutures on their boundary.
This invariant can be
seen as a natural generalization of the knot Floer homology $\HFKa$,
and its behavior under sutured decompositions can be used to prove or
re-prove important properties of knot Floer
homology~\cite{Juhasz08:SuturedDecomp}. (See
also~\cite{KronheimerMrowka10:sutured} for an analogous construction
in gauge theory.)

Ruman Zarev realized both bordered Floer homology and sutured
Floer homology as special cases of a more general theory, bordered-sutured
Floer homology~\cite{Zarev09:BorSut,Zarev:JoinGlue}.
His framework
unifies many different objects appearing in bordered Floer homology,
and the geometric objects it allows are important for proving
properties of the theory, like the formulation of the pairing theorem
in terms of morphism spaces described above, or John Etnyre, Shea
Vela-Vick, and Zarev's description of $\HFKm$ in terms of
$\HFKa$~\cite{EVVZ17} (see
also~\cite{LS:contact-glue}). Bordered-sutured Floer homology
also gives an invariant of tangles (strong enough, for instance, to
detect trivial tangles~\cite{AL19:incompressible}); see also
Section~\ref{sec:tangles}.

\subsection{Cornered Floer homology}
Since the Sei\-berg-Witten invariant was shown to admit extensions to 3\hyp manifolds
(Heegaard Floer homology, monopole Floer homology) and then 2\hyp manifolds (bordered Floer
homology), it is natural to ask if it can be extended
further. Christopher Douglas and Manolescu proposed the next step, an
extension to closed 1\hyp manifolds, surfaces with boundary, and
3\hyp manifolds with corners, which they called \emph{cornered Floer
  homology}~\cite{DouglasManolescu11:cornered}. The invariant of a
circle is a kind of algebra with both horizontal and vertical
multiplications, which they called a \emph{sequential 2\hyp algebra}. To a
surface with boundary they associate an algebra-module over this
(which behaves like an algebra with respect to horizontal
multiplication, say, and a module with respect to vertical
multiplication), and to a 3\hyp manifold with a corner, a
2\hyp module. Suitable tensor products recover bordered Floer
homology. Their construction of the 3\hyp manifold invariants uses
Sarkar-Wang's notion of nice diagrams, which makes the construction
combinatorial but means they were not able to prove invariance
directly. A variant of their construction, with somewhat more
complicated structures but where it was possible to prove invariance,
was given by Douglas, Manolescu, and the first
author~\cite{DLM:cornered}. A key idea is that the invariants of a
3\hyp manifold $Y$ with boundary $F_1\cup_{S^1}F_2$ and a corner at $S^1$
can be obtained from the bordered invariant of the smoothing of $Y$ by
taking a tensor product with an invariant of
$[0,1]\times (F_1\cup_{S^1}F_2)$, viewed as smooth on one side and
with a corner on the other. Invariance then follows from the
invariance theorem for bordered Floer homology, and the pairing
theorem follows from a computation for certain specific tri-modules and the
bordered Floer pairing theorem.

Work of Andrew Manion and Raphael Rouquier has revealed surprising
connections between the original Douglas-Manolescu construction and
representation theory~\cite{ManionRouquier}.

\subsection{Invariants of contact manifolds}
As mentioned above, a contact structure $\xi$ on a closed $3$-manifold
$Y$ induces a class $c(\xi)\in\HFa(-Y)$. It is natural to guess that if
$-Y$ is decomposed as $Y_1\cup_FY_2$ then there should be classes
$c(\xi_1)\in \CFAa(Y_1)$ and $c(\xi_2)\in\CFDa(Y_2)$ so that, under
the pairing theorem, $c(\xi)=c(\xi_1)\otimes c(\xi_2)$. Akram
Alishahi, Vikt\'{o}ria F\"{o}ldv\'{a}ri, Kristen Hendricks, Joan
Licata, Petkova, and Vera V\'{e}rtesi~\cite{AFHLPV23:contact} showed
that this is, in fact the case. Key to their definition is a precise
formulation of how $\xi$ and $F$ should interact or, equivalently,
what structure $\xi_i$ should induce on the boundary of $Y_i$. The
answer is given by the notion of a \emph{foliated open book},
introduced by Licata and V\'{e}rtesi~\cite{LicataVertesi:foliated-OB}. A
second ingredient is Honda, William Kazez, and Gordana Matic's definition of
$c(\xi)$~\cite{HKM09:contact}, which is tied more directly to a
Heegaard diagram than the original definition, and a third is the
flexibility provided by Zarev's bordered-sutured theory. A key
property of $c(\xi)$ is that $c(\xi)=0$ if $\xi$ is overtwisted (the
class of contact structures satisfying an
$h$-principle~\cite{Eliashberg89:classify-OT}). The bordered contact
invariant leads to a satisfying new proof of this fact: it reduces to
a local computation of the contact invariant near an overtwisted
disk~\cite{AHLPV:gentle-contact}.

There are also other hints of connections between bordered Floer
homology and contact topology. For instance, in unpublished work Honda
associated a triangulated category to a surface, which he called the
\emph{contact category}; objects encode contact structures near a
surface and the morphism spaces are generated by bypass
attachments. Benjamin Cooper showed that a version of the contact
category maps to the category of modules over the bordered-sutured
algebras, and at least in special cases this map is an
equivalence~\cite{Cooper:contact-cat}.  Daniel Mathews gave another
result along these lines~\cite{Mathews19:strand}; as he notes, this
implies there are $\Ainf$-style operations on the set of contact
structures, which can be computed but are not yet understood
geometrically~\cite{Mathews21:strand-Ainf}.

\subsection{Invariants of tangles}\label{sec:tangles}
Bordered Floer homology was partly inspired by extensions of Khovanov
homology to tangles (by Khovanov~\cite{Khovanov02:Tangles} and Dror
Bar-Natan~\cite{BarNatan05:Kh-tangle-cob}), and attempts to give an
extension of $\HFKa$ to tangles led to many of the basic definitions
in bordered Floer homology~\cite{LOT0}. Nonetheless, technical
difficulties relating to invariance prevented us from carrying out
this construction. Since then, three other extensions of knot Floer
homology, in the spirit of bordered Floer homology, have emerged. The
first was mentioned above: tangle exteriors are a special case of
Zarev's bordered-sutured Floer homology. A second construction was
given by Petkova-V\'ertesi~\cite{PV16:tangle-Floer}. Like our earlier attempt, they
start from a variant of Manolescu-Ozsv\'ath-Sarkar's grid diagrams. In
particular, the definition of their invariants is combinatorial and
their relationship to knot Floer homology is immediate. By contrast,
their proof of invariance uses holomorphic curves. (Indeed, for the
most elaborate version of their construction, invariance remains a
conjecture: proving it requires overcoming analytic obstacles similar
to those for bordered $\HFm$ discussed in Section~\ref{sec:Bminus}
below.) A third tangle invariant is due to Szab\'o and the second
author. It starts from a standard knot diagram, and the Heegaard
diagram it induces, and associates $\Ainf$-bimodules over certain
algebras to cups, caps, and crossings. Tensoring these bimodules together, one obtains chain complexes
associated to knot diagrams. Two papers~\cite{OSBorderedKnots1,OSBorderedKnots2}
give algebraic definitions of these bimodules and prove that the
homology groups
of the resulting chain complexes are indeed knot invariants.
In another paper~\cite{OSBorderedKnots3}, it is shown that
these operations
correspond to holomorphic curve counts and, using this, show that the
theory indeed recovers knot Floer homology. One remarkable
feature of this extension is how efficient it is for computation: while
the previous algorithm for computing $\HFKa$
using grid diagrams is effective only up to 16
crossings or so, this tangle invariant allows computations for many
knots of 80 or more crossings~\cite{Szabo:code}.  

\subsection{Fukaya-categorical invariants}\label{sec:immersed}
If $F$ is a torus $T^2$, the algebra $A(T^2)$ is an explicit,
8-dimensional algebra over $\FF_2$ with trivial differential. The
invariant of a manifold with torus boundary is a differential module
over this algebra. Jonathan Hanselman, Rasmussen, and Liam
Watson~\cite{HRW} showed that this invariant is equivalent to a
simple, geometric object: an immersed 1\hyp manifold in $T^2$, perhaps
equipped with a local system (over $\FF_2$). The tensor product in the
pairing theorem turns into taking Floer homology of curves in the
torus, an entirely combinatorial construction. The dependence of the
invariant on the parametrization of the boundary becomes transparent:
if the immersed curve associated to $Y$ is viewed as lying in
$\bdy Y$, no parametrization of the boundary is needed. (So, in a
sense, the theory becomes borderless.)

While there are reasons one might expect such a
result~\cite{LekiliPerutz11:torus,Auroux10:ICM,HKK17:immersed-curves},
their construction is complicated, surprising---and useful. In
particular, in addition to obvious computational applications, it also
has theoretical ones. One of the most remarkable is to Cameron
Gordon's Cosmetic Surgery Conjecture~\cite{CosmeticSurgeries}, which states that
distinct Dehn surgeries on a knot $K\subset S^3$ never give the same
oriented 3\hyp manifold.
Although there had been earlier progress on this problem using Heegaard Floer
techniques~\cite{OS11:RatSurg,Wang06:cosmetic,Wu11:cosmetic-L,NiWu15:cosmetic-S3,Gainullin17:mapping-cone}, Hanselman's approach
provided a powerful new framework for studying this topological problem.
Using this theory of immersed curves, Hanselman
showed that if $r$-surgery and $s$-surgery on $K$ are
orientation-preserving homeomorphic then $\{r,s\}$ is either
$\{\pm 2\}$ or $\{\pm 1/q\}$ for a particular integer $q$, as well as
a number of further restrictions (e.g., that in the first case the
knot has genus $2$)~\cite{Hanselman23:cosmetic}. Konstantinos
Varvarezos has shown that if one drops the requirement that the
homeomorphism be orientation-preserving (in which case there are many
known examples of cosmetic surgery,
such as~\cite{BHW99:chiral-cosmetic-counter}), the
immersed curve theory still
gives obstructions~\cite{Var:cosmetic}.

There are also bordered-style invariants via immersed curves in other
settings. Claudius Zibrowius introduced an immersed curve invariant
for 4\hyp ended tangles extending knot Floer homology and used it to show
that a version of knot Floer homology is invariant under
mutation~\cite{Zibrowius:mutation}, verifying a conjecture formulated
by John Baldwin and
Adam Levine~\cite{BaldwinLevine12:spanning}. Tye Lidman, Allison
Moore, and Zibrowius used this
invariant to show that so-called $L$-space knots have no essential
Conway spheres~\cite{LMZ22:no-Conway}, verifying a conjecture Lidman
and Moore formulated seven years
earlier~\cite{LidmanMoore16:pretzel}. (See~\cite{HeddenLevine16:splicing,Eftekhary18:tori,HRW}
for analogous results for closed 3\hyp manifolds.)
In a different direction, Artem
Kotelskiy, Watson, and
Zibrowius showed that, for 4\hyp ended tangles, Bar-Natan's extension of
Khovanov homology to tangles can also be interpreted as an immersed
curve in a 4\hyp punctured sphere~\cite{KLZ:immerse-Khovanov}, and this
immersed curve in fact agrees~\cite{KLZ:Kh-agree} with an invariant
introduced by Hedden, Christopher Herald, Matthew Hogancamp, and Paul
Kirk, inspired by instanton link homology~\cite{HHHK20:Fuk-SU2} (see
also~\cite{HHK14:pillow,Kotelskiy19:bord-pillow}).

\section{Bordered \texorpdfstring{$\HFm$}{HF-minus}}\label{sec:Bminus}
Extending bordered Floer homology to the $\HFm$ variant of Heegaard
Floer homology presents a set of interesting algebraic, geometric, and
analytical challenges.  To date, most work has focused on the case
where the boundary is $T^2$.

In bordered Floer homology,
the operations on the (bordered) modules are defined by counting rigid
holomorphic curves with boundary on certain Lagrangian
submanifolds. These submanifolds are non-compact; elements of the
bordered algebra record the possible asymptotics of curves, while
operations
on the algebra correspond to (most) codimension-1 degenerations of
in the moduli spaces. Passing from $\HFa$ to $\HFm$ introduces two new
complications. First, there are complicated, new codimension-1
degenerations;
philosophically, these correspond to non-empty moduli
spaces of holomorphic disks with boundary on a single Lagrangian
(cf.~\cite{FOOO1,FOOO2}). Second, there are new possible asymptotics
of curves, coming from
orbits going off to the boundary. The result is that the
algebra $A(T^2)$ must be replaced by a formal 1\hyp parameter deformation
of an $\Ainf$-algebra, or a \emph{weighted
  $\Ainf$-algebra}~\cite{LOT:torus-alg}. The higher operations
correspond to the new kinds of degenerations, and the deformation
parameter corresponds to the number of orbits. We denote
this formal deformation by $\MAlg$. Similarly, the module $\CFAa(Y_1)$
is replaced by a \emph{weighted $\Ainf$-module} $\CFAm(Y_1)$ over
$\MAlg$~\cite{LOT:torus-mod}.

The module $\CFDa(Y_2)$ can be defined as
\[
\CFDa(Y_2)=\CFAa(Y_2)\otimes\CFDDa(\Id)
\]
where $\CFDDa(\Id)$ is a type \DD\ structure (bimodule version of a
twisted complex) associated to the identity cobordism of the surface
$F$. For the case of $\HFm$, the analogous definition has a subtlety: to
formulate the notion of a type \DD\ structure over a pair of
(weighted) $\Ainf$-algebras requires a well-behaved tensor product of
(weighted) $\Ainf$-algebras, and taking one-sided tensor products
requires some further variants of this notion. In the unweighted case,
the (surprisingly subtle) construction of the tensor product of two
$A_\infty$ algebras was given by Samson Saneblidze and Ronald
Umble~\cite{SU04:Diagonals}; with some effort, their notion extends to
the versions needed for bordered Floer
theory~\cite{LOT:abstract}. With this algebra in hand, it is easy to
construct an appropriate type \DD\ bimodule $\CFDDm(\Id)$ and define
\[
\CFDm(Y_2)=\CFAm(Y_2)\otimes\CFDDm(\Id).
\]

The next step is to establish a pairing theorem for the type $A$ and
the type $D$ modules associated to bordered manifolds.  Once again,
this pairing theorem can be viewed as a deformation of the diagonal.
In the case of $\HFa$, this deformation took place in the product of
associahedra, which is a compactification of a configuration
of points in~$\RR$ relevant to $\Ainf$-algebras; in the case of
bordered $\HFm$, the associahedron is
replaced by a suitable compactification of the configuration of points
in the interior and boundary of a disk, which we call the {\em
  associaplex}~\cite{LOT:abstract}.  Details of this pairing theorem
are forthcoming~\cite{LOT:torus-pairing}.

A different construction of an extension of $\HFm$ to bordered
3\hyp manifolds with (perhaps several) torus boundary components has been
given by Ian Zemke \cite{Zemke:bordered,Zemke:lattice}, using the link
surgery formula introduced by Manolescu and the second
author~\cite{ManolescuOzsvath:surgery}. A module over the algebra
Zemke associates to a torus captures, in a succinct way, the data
needed for the surgery formula. To construct the module associated to
a 3\hyp manifold $Y$ with torus boundary, he presents $Y$ as a Dehn
filling of a link complement in $S^3$. The link surgery formula can be
applied further to prove a version of the pairing theorem. At the time
of writing, invariance of his modules remains conjectural, but he was
already able to use them to give a sought-after proof of a 
conjecture~\cite{Zemke:lattice} about the Heegaard Floer homology of
plumbed 3\hyp manifolds~\cite{SomePlumb,Nemthi08:normal}.

A third construction of a bordered extension of $\HFm$ was recently
given by Hanselman~\cite{Hanselman:HFm}, as an extension of the
immersed curve invariants discussed in Section~\ref{sec:immersed}. In
this case, the invariants take the form of immersed curves equipped
with bounding cochains (and some additional data). He proves that, up
to appropriate equivalence, this data depends only on the bordered
$3$-manifold (which he thinks of as a 3\hyp manifold with a specified knot
in it), and that it determines the $\HFm$ invariant of any Dehn
filling of the 3\hyp manifold. At the time of writing, a formula for the
behavior when one changes the parameterization of the boundary, a
general pairing theorem for this invariant, and its relationship to
the earlier immersed curve invariants, remain conjectural, but the
underlying techniques are instrumental in his results on the Cosmetic
Surgery Conjecture mentioned in Section~\ref{sec:immersed}.

\section{Applications to 4-dimensional topology}

\subsection{Knot concordance}

Heegaard Floer homology is a useful topological tool for studying the
interplay between knot theory and smooth 4\hyp dimensional topology. One
place where this interplay is particularly visible is the study of
{\em knot concordance}. We review this subject briefly; for a more
thorough overview, see~\cite{LivingstonSurvey}.

Recall that knots $K_0$ and $K_1$ are said to be {\em
  smoothly concordant} if there is an embedded annulus $A$ in
$[0,1]\times \RR^3$ which, for $i=0,1$, meets the corresponding
boundary component $\{i\}\times \RR^3$ along $K_i$.  The set of
equivalence classes of knots can be made into an abelian group, so
that addition corresponds to connected sum. This group $\mathcal{C}$
is called the {\em smooth concordance group}, and knots that are smoothly
concordant to the unknot are called {\em smoothly slice
  knots}. 

Surprisingly little is known about this concordance group. Currently,
the only known source of torsion in $\mathcal{C}$ is knots which are
the same as their mirror images (i.e., are \emph{negatively
  amphicheiral}). In particular, it is unknown whether $\mathcal{C}$
contains elements of odd order; indeed, at the time of writing, it is conceivable that
$\mathcal{C}\cong \ZZ^{\infty}\oplus (\Zmod{2})^{\infty}$.

There is a weaker notion of concordance, where the annulus $A$, rather
than being required to be smoothly embedded, is merely
topologically flat (is covered by neighborhoods homeomorphic as pairs
to $(\RR^4,\RR^2)$). There is a corresponding group $\mathcal{C}^{\mathit{top}}$,
the {\em topological concordance group}. Once again, knots that are
topologically flatly concordant to the unknot are called {\em
  topologically slice knots}.

There is a natural quotient map
$\mathcal{C}\to {\mathcal C}^{\mathit{top}}$. The kernel consists of
knots which are topologically but not smoothly slice.  A deep theorem
of Freedman~\cite{FreedmanQuinn}
states that any knot with trivial Alexander polynomial is
topologically slice. Since the advent of gauge
theory~\cite{GompfConc}, it has been known that there are knots which
are topologically but not smoothly slice; indeed, using gauge
theory, Hisaaki Endo demonstrated an infinite set of linearly
independent knots in $\mathcal{C}$ which are topologically
slice~\cite{Endo95:inf-rank-subgp}, an especially tangible
manifestation of the richness of smooth 4\hyp manifold topology.

Knot Floer homology can also be used to study concordance phenomena.
In~\cite{Hom15:inf-rank}, Jennifer Hom constructs a different infinite
set of topologically slice knots that are linearly independent in the
smooth concordance group. Her knots are built via satellite
constructions (cabling, Whitehead doubling, and connected sums).
In~\cite{OSS17:upsilon}, Andr{\'a}s Stipsicz, Szab{\'o}, and the
second author construct infinitely many homomorphisms from the smooth
concordance group to $\ZZ$, demonstrating an infinite rank free direct
summand in the concordance group of topologically slice knots; compare
also~\cite{DHST21:more-Hom}.  The homomorphisms are constructed via
knot Floer homology, and the computations of their invariants rest on
bordered techniques for computing knot Floer homology of satellite
knots.

Another exciting application of satellite computations was given by
Levine~\cite{Levine16:nonsurjective}. Using bordered Floer homology,
he showed that the Heegaard Floer $\tau$ and $\epsilon$ invariants
behave in a predictable under Mazur-pattern satellites. Using this, he
showed that most of these satellites do not bound disks in any
rational homology ball (resolving~\cite[Problem
1.45]{Kirby97:problems}), and that there are knots in homology 3-spheres
that do not bound PL disks in any homology 4-balls
(resolving~\cite[Problem 1.31]{Kirby97:problems}).

\subsection{Exotic phenomena}
Suppose that $Y_1$ and $Y_2$ are bordered 3\hyp manifolds with boundary
parameterized by $F$ and $W$ is a 4\hyp manifold with boundary identified
with $(-Y_1)\cup [0,1]\times F\cup Y_2$.
It is natural to expect that
to this data, bordered Floer homology would associate a map
\[
  F_W\co \CFDa(Y_1)\to\CFDa(Y_2)
\]
(or similarly for $\CFAa$ or the
minus variants); and that these maps would satisfy the obvious
analogue of the pairing theorem. While it is straightforward to
associate a map $F_W$ to $W$, independence of the map from the choices
in its construction has in general not been verified (but
see~\cite{TovaBrownThesis}).

Nonetheless, bordered Floer homology has been used recently for
dramatic 4\hyp dimensional applications. For instance, Gary Guth showed
that, for any $n$, there are surfaces
$(S_1,\bdy S_1)$ and $(S_2,\bdy S_2)$ in $(B^4,S^3)$, which are topologically
isotopic rel boundary but which remain smoothly non-isotopic even
after attaching $n$ $1$-handles~\cite{Guth:exotic-stab}. (Any pair of
surfaces become isotopic after attaching enough
$1$-handles~\cite{HK79:unknot-surfaces,BS16:stab}.) His construction
starts from an example of Kyle Hayden's, who gave a pair of disks in
the 4\hyp ball which are smoothly but not topologically
isotopic~\cite{Hayden:disks,Hayden:corks}, and uses a lower bound on
stabilization distance coming from Heegaard Floer homology, discovered
by Juh\'asz and Zemke~\cite{JZ:stab-dist}. Subsequently, Sungkyung Kang
gave a pair of contractible $4$-manifolds which remain
non\hyp diffeomorphic after taking the connected sum with
$S^2\times S^2$~\cite{Kang:exotic-stab}. (Several other groundbreaking
``one-is-not-enough'' results have also appeared recently, using other
techniques.)

Both of these results require complicated computations of Heegaard
Floer invariants. The disks in Guth's result are cables of Hayden's
disks. To compute their Heegaard Floer homology, he notes that if
$F\co K_1\to K_2$ is a concordance, $P$ is a pattern in the solid
torus, and $F_P$ is the satellite of $F$, then there is some map
$F_{\#}\co \CFAa(K_1)\to\CFAa(K_2)$ so that the Heegaard Floer map
$(F_P)_*\co
\HFKm((K_1)_P)\to\HFKm((K_2)_P)$~\cite{Zemke19:link-cob,AE20:tangle-cob}
is induced by
\[
  F_{\#}\otimes\Id_{\CFDa(S^1\times D^2,P)}
\]
and the pairing theorem. The key to proving this is the pairing
theorem for triangles~\cite{LOT:DCov2}. It is not needed for this
computation that $F_{\#}$ be an invariant of $F$, just that some map
with this property exists. Kang's elaborate computation uses similar
ideas, as well as work of Hendricks and the first
author~\cite{HL19:bord-involutive} and
Kang's~\cite{Kang:bord-involutive} on a bordered extension of
Hendricks-Manolescu's involutive Floer
homology~\cite{HM17:involutive}.

Both Guth's and Kang's results can be seen as resolving relative
versions of an old question of Wall's: how many stabilizations
(connected sums with $S^2\times S^2$) are required to make a pair of
homeomorphic (or homotopy equivalent), simply-connected, closed
4\hyp manifolds diffeomorphic~\cite{Wall:on-simply}.

Other applications of Heegaard Floer homology to 4\hyp manifolds seem
around the corner. For instance, Jesse Cohen recently gave an
algorithm for computing the maps on $\HFa$ associated to cobordisms
using bordered Floer homology~\cite{Cohen:composition}, extending our
algorithm for computing $\HFa$ itself~\cite{LOT4}. (The pairing
theorem for triangles~\cite{LOT:DCov2} is again a key step; another is
a rigidity result for the modules associated to
handlebodies~\cite{HL19:bord-involutive}.) As noted above, these maps
on $\HFa$
do not provide interesting invariants of closed 4\hyp manifolds, but for
4\hyp manifolds with boundary, which have seem an explosion of interest,
they do.

Most recently, Levine, Tye
Lidman, and Lisa Piccirillo have used the immersed curve formulation
of bordered Heegaard Floer homology
to give many exotic 4\hyp manifolds with $b_1=1$~\cite{LLP:exotic}. In fact, they show
that Fintushel-Stern knot surgery~\cite{FS98:knot-surg} for any nontrivial knot in
$S^3$---including knots with trivial Alexander polynomial---can be
used to produce exotic manifolds.

\section*{Acknowledgments.}
This note was written for the proceedings
of the International Congress of Basic Science, held in Beijing in
July 2023. We thank the conference organizers for their hospitality
and for the opportunity to further share this research.

\bibliographystyle{../common/hamsplain}\bibliography{../common/heegaardfloer}

\providecommand{\bysame}{\leavevmode\hbox to3em{\hrulefill}\thinspace}
\providecommand{\MR}{\relax\ifhmode\unskip\space\fi MR }
\providecommand{\MRhref}[2]{%
  \href{http://www.ams.org/mathscinet-getitem?mr=#1}{#2}
}
\providecommand{\href}[2]{#2}
\providecommand{\eprint}{\begingroup \urlstyle{rm}\Url}
\begin{thebibliography}{100}

\bibitem{AE20:tangle-cob}
Akram Alishahi and Eaman Eftekhary, \emph{Tangle {F}loer homology and
  cobordisms between tangles}, J. Topol. \textbf{13} (2020), no.~4, 1582--1657.

\bibitem{AFHLPV23:contact}
Akram Alishahi, Vikt\'{o}ria F\"{o}ldv\'{a}ri, Kristen Hendricks, Joan Licata,
  Ina Petkova, and Vera V\'{e}rtesi, \emph{Bordered {F}loer homology and
  contact structures}, Forum Math. Sigma \textbf{11} (2023), Paper No. e30, 44.

\bibitem{AHLPV:gentle-contact}
Akram Alishahi, Joan Licata, Ina Petkova, and Vera V\'{e}rtesi, \emph{A
  friendly introduction to the bordered contact invariant}, 2021,
  \eprint{arXiv:2104.07616}.

\bibitem{AL19:incompressible}
Akram Alishahi and Robert Lipshitz, \emph{Bordered {F}loer homology and
  incompressible surfaces}, Ann. Inst. Fourier (Grenoble) \textbf{69} (2019),
  no.~4, 1525--1573.

\bibitem{AndersenBenePenner09:MCGroupoid}
J{\o}rgen~Ellegaard Andersen, Alex~James Bene, and R.~C. Penner, \emph{Groupoid
  extensions of mapping class representations for bordered surfaces}, Topology
  Appl. \textbf{156} (2009), no.~17, 2713--2725.

\bibitem{Auroux10:ICM}
Denis Auroux, \emph{Fukaya categories and bordered {H}eegaard-{F}loer
  homology}, Proceedings of the {I}nternational {C}ongress of {M}athematicians.
  {V}olume {II} (New Delhi), Hindustan Book Agency, 2010, pp.~917--941.

\bibitem{Auroux10:Bordered}
\bysame, \emph{Fukaya categories of symmetric products and bordered
  {H}eegaard-{F}loer homology}, J. G\"okova Geom. Topol. GGT \textbf{4} (2010),
  1--54, \eprint{arXiv:1001.4323}.

\bibitem{BaldwinGillam12:compute}
John~A. Baldwin and William~D. Gillam, \emph{Computations of {H}eegaard-{F}loer
  knot homology}, J. Knot Theory Ramifications \textbf{21} (2012), no.~8,
  1250075, 65.

\bibitem{BaldwinLevine12:spanning}
John~A. Baldwin and Adam~Simon Levine, \emph{A combinatorial spanning tree
  model for knot {F}loer homology}, Adv. Math. \textbf{231} (2012), no.~3-4,
  1886--1939.

\bibitem{BS21:same-contact}
John~A. Baldwin and Steven Sivek, \emph{On the equivalence of contact
  invariants in sutured {F}loer homology theories}, Geom. Topol. \textbf{25}
  (2021), no.~3, 1087--1164.

\bibitem{BarNatan05:Kh-tangle-cob}
Dror Bar-Natan, \emph{Khovanov's homology for tangles and cobordisms}, Geom.
  Topol. \textbf{9} (2005), 1443--1499.

\bibitem{BS16:stab}
R.~\.{I}nan\c{c} Baykur and Nathan Sunukjian, \emph{Knotted surfaces in
  4-manifolds and stabilizations}, J. Topol. \textbf{9} (2016), no.~1,
  215--231.

\bibitem{Beliakova10:grid}
Anna Beliakova, \emph{A simplification of combinatorial link {F}loer homology},
  J. Knot Theory Ramifications \textbf{19} (2010), no.~2, 125--144.

\bibitem{Bene08:ChordDiagrams}
Alex~James Bene, \emph{A chord diagrammatic presentation of the mapping class
  group of a once bordered surface}, Geom. Dedicata \textbf{144} (2010),
  171--190, \eprint{arXiv:0802.2747}.

\bibitem{BHW99:chiral-cosmetic-counter}
Steven~A. Bleiler, Craig~D. Hodgson, and Jeffrey~R. Weeks, \emph{Cosmetic
  surgery on knots}, Proceedings of the {K}irbyfest ({B}erkeley, {CA}, 1998),
  Geom. Topol. Monogr., vol.~2, Geom. Topol. Publ., Coventry, 1999, pp.~23--34.

\bibitem{Bloom11:ss}
Jonathan~M. Bloom, \emph{A link surgery spectral sequence in monopole {F}loer
  homology}, Adv. Math. \textbf{226} (2011), no.~4, 3216--3281.

\bibitem{BEHWZ03:CompactnessInSFT}
Fr\'ed\'eric Bourgeois, Yakov Eliashberg, Helmut Hofer, Kris Wysocki, and
  Eduard Zehnder, \emph{Compactness results in symplectic field theory}, Geom.
  Topol. \textbf{7} (2003), 799--888, \eprint{arXiv:math.SG/0308183}.

\bibitem{BraamDonaldson95:triangle}
P.~J. Braam and S.~K. Donaldson, \emph{Floer's work on instanton homology,
  knots and surgery}, The {F}loer memorial volume, Progr. Math., vol. 133,
  Birkh\"{a}user, Basel, 1995, pp.~195--256.

\bibitem{TovaBrownThesis}
Tova Helen~Fell Brown, \emph{Bordered {H}eegaard {F}loer homology and
  four-manifolds with corners}, Ph.D. thesis, Massachusetts Institute of
  Technology, 2011.

\bibitem{Cohen:composition}
Jesse Cohen, \emph{Composition maps in {H}eegaard {F}loer homology}, 2023,
  \eprint{arXiv:2301.08882}.

\bibitem{ColinGhigginiHonda11:HF-ECH-3}
Vincent Colin, Paolo Ghiggini, and Ko~Honda, \emph{The equivalence of
  {H}eegaard {F}loer homology and embedded contact homology {III}: from hat to
  plus}, 2012, \eprint{arXiv:1208.1526}.

\bibitem{ColinGhigginiHonda11:HF-ECH-1}
\bysame, \emph{The equivalence of {H}eegaard {F}loer homology and embedded
  contact homology via open book decompositions {I}}, 2012,
  \eprint{arXiv:1208.1074}.

\bibitem{ColinGhigginiHonda11:HF-ECH-2}
\bysame, \emph{The equivalence of {H}eegaard {F}loer homology and embedded
  contact homology via open book decompositions {II}}, 2012,
  \eprint{arXiv:1208.1077}.

\bibitem{Cooper:contact-cat}
Benjamin Cooper, \emph{Formal contact categories}, 2015,
  \eprint{arXiv:1511.04765}.

\bibitem{DHST21:more-Hom}
Irving Dai, Jennifer Hom, Matthew Stoffregen, and Linh Truong, \emph{More
  concordance homomorphisms from knot {F}loer homology}, Geom. Topol.
  \textbf{25} (2021), no.~1, 275--338.

\bibitem{DonaldsonPolynomials}
S.~K. Donaldson, \emph{Polynomial invariants for smooth four-manifolds},
  Topology \textbf{29} (1990), no.~3, 257--315.

\bibitem{Donaldson96:survey}
\bysame, \emph{The {S}eiberg-{W}itten equations and {$4$}-manifold topology},
  Bull. Amer. Math. Soc. (N.S.) \textbf{33} (1996), no.~1, 45--70.

\bibitem{DonaldsonFloer}
\bysame, \emph{Floer homology groups in {Y}ang-{M}ills theory}, Cambridge
  Tracts in Mathematics, vol. 147, Cambridge University Press, Cambridge, 2002,
  With the assistance of M. Furuta and D. Kotschick.

\bibitem{DLM:cornered}
Christopher~L. Douglas, Robert Lipshitz, and Ciprian Manolescu, \emph{Cornered
  {H}eegaard {F}loer homology}, Mem. Amer. Math. Soc. \textbf{262} (2019),
  no.~1266, v+124.

\bibitem{DouglasManolescu11:cornered}
Christopher~L. Douglas and Ciprian Manolescu, \emph{On the algebra of cornered
  {F}loer homology}, J. Topol. \textbf{7} (2014), no.~1, 1--68.

\bibitem{Dowlin:kh-to-hfk}
Nathan Dowlin, \emph{A spectral sequence from {K}hovanov homology to knot
  {F}loer homology}, 2018, \eprint{arXiv:1811.07848}.

\bibitem{Eftekhary05:LongitudeWhitehead}
Eaman Eftekhary, \emph{Longitude {F}loer homology and the {W}hitehead double},
  Algebr. Geom. Topol. \textbf{5} (2005), 1389--1418,
  \eprint{arXiv:math.GT/0407211}.

\bibitem{Eftekhary18:tori}
\bysame, \emph{Bordered {F}loer homology and existence of incompressible tori
  in homology spheres}, Compos. Math. \textbf{154} (2018), no.~6, 1222--1268.

\bibitem{Eliashberg89:classify-OT}
Y.~Eliashberg, \emph{Classification of overtwisted contact structures on
  {$3$}-manifolds}, Invent. Math. \textbf{98} (1989), no.~3, 623--637.

\bibitem{EGH00:IntroductionSFT}
Yakov Eliashberg, Alexander Givental, and Helmut Hofer, \emph{Introduction to
  symplectic field theory}, Geom. Funct. Anal. (2000), no.~Special Volume, Part
  II, 560--673, \eprint{arXiv:math.SG/0010059}, GAFA 2000 (Tel Aviv, 1999).

\bibitem{ET98:confoliations}
Yakov~M. Eliashberg and William~P. Thurston, \emph{Confoliations}, University
  Lecture Series, vol.~13, American Mathematical Society, Providence, RI, 1998.

\bibitem{Endo95:inf-rank-subgp}
Hisaaki Endo, \emph{Linear independence of topologically slice knots in the
  smooth cobordism group}, Topology Appl. \textbf{63} (1995), no.~3, 257--262.

\bibitem{EVVZ17}
John~B. Etnyre, David~Shea Vela-Vick, and Rumen Zarev, \emph{Sutured {F}loer
  homology and invariants of {L}egendrian and transverse knots}, Geom. Topol.
  \textbf{21} (2017), no.~3, 1469--1582.

\bibitem{FS98:knot-surg}
Ronald Fintushel and Ronald~J. Stern, \emph{Knots, links, and {$4$}-manifolds},
  Invent. Math. \textbf{134} (1998), no.~2, 363--400.

\bibitem{Floer88:LagrangianHF}
Andreas Floer, \emph{Morse theory for {L}agrangian intersections}, J.
  Differential Geom. \textbf{28} (1988), no.~3, 513--547.

\bibitem{FloerTriangles}
Andreas Floer, \emph{Instanton homology and {D}ehn surgery}, The {F}loer
  memorial volume, Progr. Math., no. 133, Birkh{\"a}user, 1995, pp.~77--97.

\bibitem{FreedmanQuinn}
Michael~H. Freedman and Frank Quinn, \emph{Topology of 4-manifolds}, Princeton
  Mathematical Series, vol.~39, Princeton University Press, Princeton, NJ,
  1990.

\bibitem{FOOO1}
Kenji Fukaya, Yong-Geun Oh, Hiroshi Ohta, and Kaoru Ono, \emph{Lagrangian
  intersection {F}loer theory: anomaly and obstruction. {P}art {I}}, AMS/IP
  Studies in Advanced Mathematics, vol.~46, American Mathematical Society,
  Providence, RI; International Press, Somerville, MA, 2009.

\bibitem{FOOO2}
\bysame, \emph{Lagrangian intersection {F}loer theory: anomaly and obstruction.
  {P}art {II}}, AMS/IP Studies in Advanced Mathematics, vol.~46, American
  Mathematical Society, Providence, RI; International Press, Somerville, MA,
  2009.

\bibitem{Gabai83:foliations}
David Gabai, \emph{Foliations and the topology of {$3$}-manifolds}, J.
  Differential Geom. \textbf{18} (1983), no.~3, 445--503.

\bibitem{Gabai84:knot-foliations}
\bysame, \emph{Foliations and genera of links}, Topology \textbf{23} (1984),
  no.~4, 381--394.

\bibitem{Gainullin17:mapping-cone}
Fyodor Gainullin, \emph{The mapping cone formula in {H}eegaard {F}loer homology
  and {D}ehn surgery on knots in {$S^3$}}, Algebr. Geom. Topol. \textbf{17}
  (2017), no.~4, 1917--1951.

\bibitem{Ghiggini08:FiberedGenusOne}
Paolo Ghiggini, \emph{Knot {F}loer homology detects genus-one fibred knots},
  Amer. J. Math. \textbf{130} (2008), no.~5, 1151--1169.

\bibitem{Giroux02:correspondence-ICM}
Emmanuel Giroux, \emph{G\'{e}om\'{e}trie de contact: de la dimension trois vers
  les dimensions sup\'{e}rieures}, Proceedings of the {I}nternational
  {C}ongress of {M}athematicians, {V}ol. {II} ({B}eijing, 2002), Higher Ed.
  Press, Beijing, 2002, pp.~405--414.

\bibitem{GompfConc}
Robert~E. Gompf, \emph{Smooth concordance of topologically slice knots},
  Topology \textbf{25} (1986), no.~3, 353--373.

\bibitem{CosmeticSurgeries}
Cameron~McA. Gordon, \emph{Dehn surgery on knots}, Proceedings of the
  {I}nternational {C}ongress of {M}athematicians, {V}ol. {I}, {II} ({K}yoto,
  1990), Math. Soc. Japan, Tokyo, 1991, pp.~631--642.

\bibitem{GrigsbyWehrli:detects}
J.~Elisenda Grigsby and Stephan~M. Wehrli, \emph{On the colored {J}ones
  polynomial, sutured {F}loer homology, and knot {F}loer homology}, Adv. Math.
  \textbf{223} (2010), no.~6, 2114--2165, \eprint{arXiv:0907.4375}.

\bibitem{Guth:exotic-stab}
Gary Guth, \emph{For exotic surfaces with boundary, one stabilization is not
  enough}, 2022, \eprint{arXiv:2207.11847}.

\bibitem{HKK17:immersed-curves}
F.~Haiden, L.~Katzarkov, and M.~Kontsevich, \emph{Flat surfaces and stability
  structures}, Publ. Math. Inst. Hautes \'{E}tudes Sci. \textbf{126} (2017),
  247--318.

\bibitem{Hanselman23:cosmetic}
Jonathan Hanselman, \emph{Heegaard {F}loer homology and cosmetic surgeries in
  {$S^3$}}, J. Eur. Math. Soc. (JEMS) \textbf{25} (2023), no.~5, 1627--1670.

\bibitem{Hanselman:HFm}
Jonathan Hanselman, \emph{Knot {F}loer homology as immersed curves}, 2023,
  \eprint{arXiv:2305.16271}.

\bibitem{HRW}
Jonathan Hanselman, Jacob Rasmussen, and Liam Watson, \emph{Bordered {F}loer
  homology for manifolds with torus boundary via immersed curves}, 2016,
  \eprint{arXiv:1604.03466}.

\bibitem{Hayden:disks}
Kyle Hayden, \emph{Exotically knotted disks and complex curves}, 2020,
  \eprint{arXiv:2003.13681}.

\bibitem{Hayden:corks}
\bysame, \emph{Corks, covers, and complex curves}, 2021,
  \eprint{arXiv:2107.06856}.

\bibitem{Hedden}
Matthew Hedden, \emph{On knot {F}loer homology and cabling}, Algebr. Geom.
  Topol. \textbf{5} (2005), 1197--1222, \eprint{arXiv:math.GT/0406402}.

\bibitem{HeddenWhitehead}
\bysame, \emph{Knot {F}loer homology of {W}hitehead doubles}, Geom. Topol.
  \textbf{11} (2007), 2277--2338, \eprint{arXiv:math.GT/0606094}.

\bibitem{HHHK20:Fuk-SU2}
Matthew Hedden, Christopher~M. Herald, Matthew Hogancamp, and Paul Kirk,
  \emph{The {F}ukaya category of the pillowcase, traceless character varieties,
  and {K}hovanov cohomology}, Trans. Amer. Math. Soc. \textbf{373} (2020),
  no.~12, 8391--8437.

\bibitem{HHK14:pillow}
Matthew Hedden, Christopher~M. Herald, and Paul Kirk, \emph{The pillowcase and
  perturbations of traceless representations of knot groups}, Geom. Topol.
  \textbf{18} (2014), no.~1, 211--287.

\bibitem{HeddenLevine16:splicing}
Matthew Hedden and Adam~Simon Levine, \emph{Splicing knot complements and
  bordered {F}loer homology}, J. Reine Angew. Math. \textbf{720} (2016),
  129--154, \eprint{arXiv:1210.7055}.

\bibitem{HL19:bord-involutive}
Kristen Hendricks and Robert Lipshitz, \emph{Involutive bordered {F}loer
  homology}, Trans. Amer. Math. Soc. \textbf{372} (2019), no.~1, 389--424.

\bibitem{HM17:involutive}
Kristen Hendricks and Ciprian Manolescu, \emph{Involutive {H}eegaard {F}loer
  homology}, Duke Math. J. \textbf{166} (2017), no.~7, 1211--1299.

\bibitem{Hockenhull:bord-link}
Thomas Hockenhull, \emph{Holomorphic polygons and the bordered {H}eegaard
  {F}loer homology of link complements}, 2018, \eprint{arXiv:1802.02443}.

\bibitem{Hom15:inf-rank}
Jennifer Hom, \emph{An infinite-rank summand of topologically slice knots},
  Geom. Topol. \textbf{19} (2015), no.~2, 1063--1110.

\bibitem{HKM09:contact}
Ko~Honda, William~H. Kazez, and Gordana Mati{\'c}, \emph{On the contact class
  in {H}eegaard {F}loer homology}, J. Differential Geom. \textbf{83} (2009),
  no.~2, 289--311.

\bibitem{HK79:unknot-surfaces}
Fujitsugu Hosokawa and Akio Kawauchi, \emph{Proposals for unknotted surfaces in
  four-spaces}, Osaka Math. J. \textbf{16} (1979), no.~1, 233--248.

\bibitem{Hutchings02:ECH}
Michael Hutchings, \emph{An index inequality for embedded pseudoholomorphic
  curves in symplectizations}, J. Eur. Math. Soc. (JEMS) \textbf{4} (2002),
  no.~4, 313--361.

\bibitem{HutchingsTaubes07:GluingI}
Michael Hutchings and Clifford~Henry Taubes, \emph{Gluing pseudoholomorphic
  curves along branched covered cylinders. {I}}, J. Symplectic Geom. \textbf{5}
  (2007), no.~1, 43--137.

\bibitem{HutchingsTaubes09:GluingII}
\bysame, \emph{Gluing pseudoholomorphic curves along branched covered
  cylinders. {II}}, J. Symplectic Geom. \textbf{7} (2009), no.~1, 29--133.

\bibitem{Juhasz06:Sutured}
Andr{\'a}s Juh{\'a}sz, \emph{Holomorphic discs and sutured manifolds}, Algebr.
  Geom. Topol. \textbf{6} (2006), 1429--1457, \eprint{arXiv:math/0601443}.

\bibitem{Juhasz08:SuturedDecomp}
\bysame, \emph{Floer homology and surface decompositions}, Geom. Topol.
  \textbf{12} (2008), no.~1, 299--350, \eprint{arXiv:math/0609779}.

\bibitem{JZ:stab-dist}
Andr{\'a}s Juh{\'a}sz and Ian Zemke, \emph{Stabilization distance bounds from
  link {F}loer homology}, 2018, \eprint{arXiv:1810.09158}.

\bibitem{Kang:bord-involutive}
Sungkyung Kang, \emph{Involutive knot {F}loer homology and bordered modules},
  2022, \eprint{arXiv:2202.12500}.

\bibitem{Kang:exotic-stab}
\bysame, \emph{One stabilization is not enough for contractible 4-manifolds},
  2022, \eprint{arXiv:2210.07510}.

\bibitem{Khovanov00:CatJones}
Mikhail Khovanov, \emph{A categorification of the {J}ones polynomial}, Duke
  Math. J. \textbf{101} (2000), no.~3, 359--426.

\bibitem{Khovanov02:Tangles}
\bysame, \emph{A functor-valued invariant of tangles}, Algebr. Geom. Topol.
  \textbf{2} (2002), 665--741.

\bibitem{Kirby97:problems}
\emph{Problems in low-dimensional topology}, Geometric topology ({A}thens,
  {GA}, 1993) (Rob Kirby, ed.), AMS/IP Stud. Adv. Math., vol. 2.2, Amer. Math.
  Soc., Providence, RI, 1997, pp.~35--473.

\bibitem{Kotelskiy19:bord-pillow}
Artem Kotelskiy, \emph{Bordered theory for pillowcase homology}, Math. Res.
  Lett. \textbf{26} (2019), no.~5, 1467--1516.

\bibitem{KLZ:immerse-Khovanov}
Artem Kotelskiy, Liam Watson, and Claudius Zibrowius, \emph{Immersed curves in
  {K}hovanov homology}, arXiv:1910.14584.

\bibitem{KLZ:Kh-agree}
\bysame, \emph{{K}hovanov invariants via {F}ukaya categories: the tangle
  invariants agree}, arXiv:2004.01619.

\bibitem{ManolescuKronheimerSWspectrum}
Peter Kronheimer and Ciprian Manolescu, \emph{Periodic {F}loer pro-spectra from
  the {S}eiberg-{W}itten equations}, arXiv:math/0203243, 2002.

\bibitem{KM97:contact}
Peter Kronheimer and Tomasz Mrowka, \emph{Monopoles and contact structures},
  Invent. Math. \textbf{130} (1997), no.~2, 209--255.

\bibitem{KM04:P}
\bysame, \emph{Witten's conjecture and property {P}}, Geom. Topol. \textbf{8}
  (2004), 295--310.

\bibitem{KronheimerMrowka}
\bysame, \emph{Monopoles and three-manifolds}, New Mathematical Monographs,
  vol.~10, Cambridge University Press, Cambridge, 2007.

\bibitem{KronheimerMrowka10:sutured}
\bysame, \emph{Knots, sutures, and excision}, J. Differential Geom. \textbf{84}
  (2010), no.~2, 301--364, \eprint{arXiv:0807.4891}.

\bibitem{KronheimerMrowka11:detect}
\bysame, \emph{Khovanov homology is an unknot-detector}, Publ. Math. Inst.
  Hautes \'Etudes Sci. (2011), no.~113, 97--208, \eprint{arXiv:1005.4346}.

\bibitem{KMOS}
Peter Kronheimer, Tomasz Mrowka, Peter~S. Ozsv{\'a}th, and Zolt{\'a}n {\relax
  Sz}ab{\'o}, \emph{Monopoles and lens space surgeries}, Ann. of Math. (2)
  \textbf{165} (2007), no.~2, 457--546, \eprint{arXiv:math.GT/0310164}.

\bibitem{KutluhanLeeTaubes20:HFHMIV}
\c{C}a\u{g}atay Kutluhan, Yi-Jen Lee, and Clifford~Henry Taubes,
  \emph{H{F}={HM}, {IV}: {T}he {S}ieberg-{W}itten {F}loer homology and ech
  correspondence}, Geom. Topol. \textbf{24} (2020), no.~7, 3219--3469.

\bibitem{KutluhanLeeTaubes20:HFHMV}
\bysame, \emph{H{F}={HM}, {V}: {S}eiberg-{W}itten {F}loer homology and handle
  additions}, Geom. Topol. \textbf{24} (2020), no.~7, 3471--3748.

\bibitem{KutluhanLeeTaubes20:HFHMIII}
\bysame, \emph{{$\rm HF{=}HM$}, {III}: holomorphic curves and the differential
  for the ech/{H}eegaard {F}loer correspondence}, Geom. Topol. \textbf{24}
  (2020), no.~6, 3013--3218.

\bibitem{KutluhanLeeTaubes20:HFHMI}
\bysame, \emph{{${\rm HF}={\rm HM}$}, {I}: {H}eegaard {F}loer homology and
  {S}eiberg-{W}itten {F}loer homology}, Geom. Topol. \textbf{24} (2020), no.~6,
  2829--2854.

\bibitem{KutluhanLeeTaubes20:HFHMII}
\bysame, \emph{{${\rm HF}={\rm HM}$}, {II}: {R}eeb orbits and holomorphic
  curves for the ech/{H}eegaard {F}loer correspondence}, Geom. Topol.
  \textbf{24} (2020), no.~6, 2855--3012.

\bibitem{LS:contact-glue}
Ryan Leigon and Federico Salmoiraghi, \emph{Equivalence of contact gluing maps
  in sutured {F}loer homology}, 2020, \eprint{arXiv:2005.04827}.

\bibitem{LekiliPerutz11:torus}
Yank{\i } Lekili and Timothy Perutz, \emph{Fukaya categories of the torus and
  {D}ehn surgery}, Proc. Natl. Acad. Sci. USA \textbf{108} (2011), no.~20,
  8106--8113.

\bibitem{Levine16:nonsurjective}
Adam~Simon Levine, \emph{Nonsurjective satellite operators and piecewise-linear
  concordance}, Forum Math. Sigma \textbf{4} (2016), Paper No. e34, 47.

\bibitem{LLP:exotic}
Adam~Simon Levine, Tye Lidman, and Lisa Piccirillo, \emph{New constructions and
  invariants of closed exotic 4-manifolds}, 2023, \eprint{arXiv:2307.08130}.

\bibitem{LicataVertesi:foliated-OB}
Joan Licata and Vera V\'{e}rtesi, \emph{Foliated open books}, 2020,
  \eprint{arXiv:2002.01752}.

\bibitem{LidmanMoore16:pretzel}
Tye Lidman and Allison~H. Moore, \emph{Pretzel knots with {$L$}-space
  surgeries}, Michigan Math. J. \textbf{65} (2016), no.~1, 105--130.

\bibitem{LMZ22:no-Conway}
Tye Lidman, Allison~H. Moore, and Claudius Zibrowius, \emph{{$L$}-space knots
  have no essential {C}onway spheres}, Geom. Topol. \textbf{26} (2022), no.~5,
  2065--2102.

\bibitem{LOT:torus-pairing}
Robert Lipshitz, Peter~S. Ozsv{\'a}th, and Dylan~P. Thurston, \emph{Pairing
  theorems for bordered {$\mathit{HF}^-$} with torus boundary}, In preparation.

\bibitem{LOT0}
\bysame, \emph{Slicing planar grid diagrams: a gentle introduction to bordered
  {H}eegaard {F}loer homology}, Proceedings of {G}\"okova {G}eometry-{T}opology
  {C}onference 2008, G\"okova Geometry/Topology Conference (GGT), G\"okova,
  2009, \eprint{arXiv:0810.0695}, pp.~91--119.

\bibitem{LOTHomPair}
\bysame, \emph{Heegaard {F}loer homology as morphism spaces}, Quantum Topol.
  \textbf{2} (2011), no.~4, 381--449, \eprint{arXiv:1005.1248}.

\bibitem{LOT:DCov1}
\bysame, \emph{Bordered {F}loer homology and the spectral sequence of a
  branched double cover {I}}, J. Topol. \textbf{7} (2014), no.~4, 1155--1199,
  \eprint{arXiv:1011.0499}.

\bibitem{LOT4}
\bysame, \emph{Computing {$\widehat{\mathit{HF}}$} by factoring mapping
  classes}, Geom. Topol. \textbf{18} (2014), no.~5, 2547--2681,
  \eprint{arXiv:1010.2550}.

\bibitem{LOT2}
\bysame, \emph{Bimodules in bordered {H}eegaard {F}loer homology}, Geom. Topol.
  \textbf{19} (2015), no.~2, 525--724, \eprint{arXiv:1003.0598}.

\bibitem{LOT:DCov2}
\bysame, \emph{Bordered {F}loer homology and the spectral sequence of a
  branched double cover {II}: the spectral sequences agree}, J. Topol.
  \textbf{9} (2016), no.~2, 607--686, \eprint{arXiv:1404.2894}.

\bibitem{LOT1}
\bysame, \emph{Bordered {H}eegaard {F}loer homology}, Mem. Amer. Math. Soc.
  \textbf{254} (2018), no.~1216, viii+279, \eprint{arXiv:0810.0687}.

\bibitem{LOT:abstract}
\bysame, \emph{Diagonals and {$A$}-infinity tensor products}, 2020,
  \eprint{arXiv:2009.05222}.

\bibitem{LOT:torus-alg}
\bysame, \emph{A bordered {$\mathit{HF}^-$} algebra for the torus}, 2021,
  \eprint{arXiv:2108.12488}.

\bibitem{LOT:torus-mod}
\bysame, \emph{Bordered {$\HFm$} for three-manifolds with torus boundary},
  2023, \eprint{arXiv:2305.07754}.

\bibitem{LScontactInvI}
Paolo Lisca and Andr\'{a}s~I. Stipsicz, \emph{Ozsv\'{a}th-{S}zab\'{o}
  invariants and tight contact three-manifolds. {I}}, Geom. Topol. \textbf{8}
  (2004), 925--945.

\bibitem{LScontactInv0}
\bysame, \emph{Seifert fibered contact three-manifolds via surgery}, Algebr.
  Geom. Topol. \textbf{4} (2004), 199--217.

\bibitem{LScontactInvII}
\bysame, \emph{Ozsv\'{a}th-{S}zab\'{o} invariants and tight contact
  three-manifolds. {II}}, J. Differential Geom. \textbf{75} (2007), no.~1,
  109--141.

\bibitem{LivingstonSurvey}
Charles Livingston, \emph{A survey of classical knot concordance}, Handbook of
  knot theory, Elsevier B. V., Amsterdam, 2005, pp.~319--347.

\bibitem{ManionRouquier}
Andrew Manion and Raphael Rouquier, \emph{Higher representations and cornered
  {H}eegaard {F}loer homology}, 2020, \eprint{arXiv:2009.09627}.

\bibitem{Manolescu03:SW-spectrum-1}
Ciprian Manolescu, \emph{Seiberg-{W}itten-{F}loer stable homotopy type of
  three-manifolds with {$b_1=0$}}, Geom. Topol. \textbf{7} (2003), 889--932
  (electronic).

\bibitem{ManolescuOzsvath:surgery}
Ciprian Manolescu and Peter Ozsv{\'a}th, \emph{Heegaard {F}loer homology and
  integer surgeries on links}, 2010, \eprint{arXiv:1011.1317}.

\bibitem{MOS06:CombinatorialDescrip}
Ciprian Manolescu, Peter Ozsv{\'a}th, and Sucharit Sarkar, \emph{A
  combinatorial description of knot {F}loer homology}, Ann. of Math. (2)
  \textbf{169} (2009), no.~2, 633--660, \eprint{arXiv:math.GT/0607691}.

\bibitem{MOT:grid}
Ciprian Manolescu, Peter Ozsv{\'a}th, and Dylan Thurston, \emph{Grid diagrams
  and {H}eegaard {F}loer invariants}, 2009, \eprint{arXiv:0910.0078}.

\bibitem{MOST07:CombinatorialLink}
Ciprian Manolescu, Peter~S. Ozsv{\'a}th, Zolt{\'a}n {\relax Sz}ab{\'o}, and
  Dylan~P. Thurston, \emph{On combinatorial link {F}loer homology}, Geom.
  Topol. \textbf{11} (2007), 2339--2412, \eprint{arXiv:math.GT/0610559}.

\bibitem{Mathews19:strand}
Daniel~V. Mathews, \emph{Strand algebras and contact categories}, Geom. Topol.
  \textbf{23} (2019), no.~2, 637--683.

\bibitem{Mathews21:strand-Ainf}
\bysame, \emph{A-infinity algebras, strand algebras, and contact categories},
  Algebr. Geom. Topol. \textbf{21} (2021), no.~3, 1093--1207.

\bibitem{Nemthi08:normal}
Andr\'{a}s N\'{e}methi, \emph{Lattice cohomology of normal surface
  singularities}, Publ. Res. Inst. Math. Sci. \textbf{44} (2008), no.~2,
  507--543.

\bibitem{Ni07:FiberedKnot}
Yi~Ni, \emph{Knot {F}loer homology detects fibred knots}, Invent. Math.
  \textbf{170} (2007), no.~3, 577--608.

\bibitem{NiWu15:cosmetic-S3}
Yi~Ni and Zhongtao Wu, \emph{Cosmetic surgeries on knots in {$S^3$}}, J. Reine
  Angew. Math. \textbf{706} (2015), 1--17.

\bibitem{OSS15:grid-book}
Peter~S. Ozsv\'{a}th, Andr\'{a}s~I. Stipsicz, and Zolt\'{a}n Szab\'{o},
  \emph{Grid homology for knots and links}, Mathematical Surveys and
  Monographs, vol. 208, American Mathematical Society, Providence, RI, 2015.

\bibitem{OSS17:upsilon}
\bysame, \emph{Concordance homomorphisms from knot {F}loer homology}, Adv.
  Math. \textbf{315} (2017), 366--426.

\bibitem{SomePlumb}
Peter~S. Ozsv{\'a}th and Zolt{\'a}n {\relax Sz}ab{\'o}, \emph{On the {F}loer
  homology of plumbed three-manifolds}, Geom. Topol. \textbf{7} (2003),
  185--224, \eprint{arXiv:math.SG/0203265}.

\bibitem{OS04:ThurstonNorm}
\bysame, \emph{Holomorphic disks and genus bounds}, Geom. Topol. \textbf{8}
  (2004), 311--334.

\bibitem{OS04:HolDiskProperties}
\bysame, \emph{Holomorphic disks and three-manifold invariants: properties and
  applications}, Ann. of Math. (2) \textbf{159} (2004), no.~3, 1159--1245,
  \eprint{arXiv:math.SG/0105202}.

\bibitem{OS04:HolomorphicDisks}
\bysame, \emph{Holomorphic disks and topological invariants for closed
  three-manifolds}, Ann. of Math. (2) \textbf{159} (2004), no.~3, 1027--1158,
  \eprint{arXiv:math.SG/0101206}.

\bibitem{OS05:Contact}
\bysame, \emph{Heegaard {F}loer homology and contact structures}, Duke Math. J.
  \textbf{129} (2005), no.~1, 39--61.

\bibitem{BrDCov}
\bysame, \emph{On the {H}eegaard {F}loer homology of branched double-covers},
  Adv. Math. \textbf{194} (2005), no.~1, 1--33, \eprint{arXiv:math.GT/0309170}.

\bibitem{OS06:HolDiskFour}
\bysame, \emph{Holomorphic triangles and invariants for smooth four-manifolds},
  Adv. Math. \textbf{202} (2006), no.~2, 326--400,
  \eprint{arXiv:math.SG/0110169}.

\bibitem{OS11:RatSurg}
\bysame, \emph{Knot {F}loer homology and rational surgeries}, Algebr. Geom.
  Topol. \textbf{11} (2011), no.~1, 1--68.

\bibitem{OSBorderedKnots2}
\bysame, \emph{Kauffman states, bordered algebras, and a bigraded knot
  invariant}, Adv. Math. \textbf{328} (2018), 1088--1198.

\bibitem{OSBorderedKnots3}
\bysame, \emph{Algebras with matchings and {K}not {F}loer homology}, 2019,
  \eprint{arXiv:1912.01657}.

\bibitem{OSBorderedKnots1}
\bysame, \emph{Bordered knot algebras with matchings}, Quantum Topol.
  \textbf{10} (2019), no.~3, 481--592.

\bibitem{Petkova18:decat}
Ina Petkova, \emph{The decategorification of bordered {H}eegaard {F}loer
  homology}, J. Symplectic Geom. \textbf{16} (2018), no.~1, 227--277.

\bibitem{PV16:tangle-Floer}
Ina Petkova and Vera V\'{e}rtesi, \emph{Combinatorial tangle {F}loer homology},
  Geom. Topol. \textbf{20} (2016), no.~6, 3219--3332.

\bibitem{SU04:Diagonals}
Samson Saneblidze and Ronald Umble, \emph{Diagonals on the permutahedra,
  multiplihedra and associahedra}, Homology Homotopy Appl. \textbf{6} (2004),
  no.~1, 363--411, \eprint{arXiv:math.AT/0209109}.

\bibitem{SarkarWang07:ComputingHFhat}
Sucharit Sarkar and Jiajun Wang, \emph{An algorithm for computing some
  {H}eegaard {F}loer homologies}, Ann. of Math. (2) \textbf{171} (2010), no.~2,
  1213--1236, \eprint{arXiv:math/0607777}.

\bibitem{Stasheff63:associahedron1}
James~Dillon Stasheff, \emph{Homotopy associativity of {$H$}-spaces. {I}},
  Trans. Amer. Math. Soc. \textbf{108} (1963), 275--292.

\bibitem{Szabo:code}
Zolt\'{a}n Szab\'{o}, \emph{Knot {F}loer homology calculator},
  \url{https://web.math.princeton.edu/~szabo/HFKcalc.html}.

\bibitem{Szabo15:ss}
\bysame, \emph{A geometric spectral sequence in {K}hovanov homology}, J. Topol.
  \textbf{8} (2015), no.~4, 1017--1044.

\bibitem{Taubes94:SW-symplectic}
Clifford~Henry Taubes, \emph{The {S}eiberg-{W}itten invariants and symplectic
  forms}, Math. Res. Lett. \textbf{1} (1994), no.~6, 809--822.

\bibitem{Taubes10:ECH-SW1}
\bysame, \emph{Embedded contact homology and {S}eiberg-{W}itten {F}loer
  cohomology {I}}, Geom. Topol. \textbf{14} (2010), no.~5, 2497--2581.

\bibitem{Taubes10:ECH-SW2}
\bysame, \emph{Embedded contact homology and {S}eiberg-{W}itten {F}loer
  cohomology {II}}, Geom. Topol. \textbf{14} (2010), no.~5, 2583--2720.

\bibitem{Taubes10:ECH-SW3}
\bysame, \emph{Embedded contact homology and {S}eiberg-{W}itten {F}loer
  cohomology {III}}, Geom. Topol. \textbf{14} (2010), no.~5, 2721--2817.

\bibitem{Taubes10:ECH-SW4}
\bysame, \emph{Embedded contact homology and {S}eiberg-{W}itten {F}loer
  cohomology {IV}}, Geom. Topol. \textbf{14} (2010), no.~5, 2819--2960.

\bibitem{Taubes10:ECH-SW5}
\bysame, \emph{Embedded contact homology and {S}eiberg-{W}itten {F}loer
  cohomology {V}}, Geom. Topol. \textbf{14} (2010), no.~5, 2961--3000.

\bibitem{Var:cosmetic}
Konstantinos Varvarezos, \emph{Heegaard {F}loer homology and chirally cosmetic
  surgeries}, arXiv:2112.03144.

\bibitem{Wall:on-simply}
C.~T.~C. Wall, \emph{On simply-connected {$4$}-manifolds}, J. London Math. Soc.
  \textbf{39} (1964), 141--149.

\bibitem{Wang06:cosmetic}
Jiajun Wang, \emph{Cosmetic surgeries on genus one knots}, Algebr. Geom. Topol.
  \textbf{6} (2006), 1491--1517.

\bibitem{Witten}
Edward Witten, \emph{Monopoles and four-manifolds}, Math. Res. Lett. \textbf{1}
  (1994), 769--796, \eprint{arXiv:hep-th/9411102}.

\bibitem{Wu11:cosmetic-L}
Zhongtao Wu, \emph{Cosmetic surgery in {L}-space homology spheres}, Geom.
  Topol. \textbf{15} (2011), no.~2, 1157--1168.

\bibitem{Zarev09:BorSut}
Rumen Zarev, \emph{Bordered {F}loer homology for sutured manifolds}, 2009,
  \eprint{arXiv:0908.1106}.

\bibitem{Zarev:JoinGlue}
\bysame, \emph{Joining and gluing sutured {F}loer homology}, 2010,
  \eprint{arXiv:1010.3496}.

\bibitem{Zemke:bordered}
Ian Zemke, \emph{Bordered manifolds with torus boundary and the link surgery
  formula}, arXiv:2109.11520.

\bibitem{Zemke:lattice}
\bysame, \emph{The equivalence of lattice and {H}eegaard {F}loer homology},
  arXiv:2111.14962.

\bibitem{Zemke19:link-cob}
Ian Zemke, \emph{Link cobordisms and functoriality in link {F}loer homology},
  J. Topol. \textbf{12} (2019), no.~1, 94--220.

\bibitem{Zhan:code}
Bohua Zhan, \emph{bfh{\underline{\ }}python},
  \url{https://github.com/bzhan/bfh_python}.

\bibitem{Zhan14:thesis}
\bysame, \emph{Combinatorial methods in bordered {H}eegaard {F}loer homology},
  ProQuest LLC, Ann Arbor, MI, 2014, Thesis (Ph.D.)--Princeton University.

\bibitem{Zhan16:proofs}
\bysame, \emph{Combinatorial proofs in bordered {H}eegaard {F}loer homology},
  Algebr. Geom. Topol. \textbf{16} (2016), no.~5, 2571--2636.

\bibitem{Zibrowius:mutation}
Claudius Zibrowius, \emph{On symmetries of peculiar modules; or,
  {$\delta$}-graded link {F}loer homology is mutation invariant},
  arXiv:1909.04267.

\end{thebibliography}
\end{document}